\documentclass[reqno,centertags, 12pt]{amsart}
\usepackage{amsmath,amsthm,amscd,amssymb}
\usepackage{latexsym}
-\marginparwidth \oddsidemargin -4mm \evensidemargin
\oddsidemargin

\newcommand{\bbN}{{\mathbb{N}}}
\newcommand{\bbR}{{\mathbb{R}}}
\newcommand{\bbD}{{\mathbb{D}}}

\newcommand{\bbZ}{{\mathbb{Z}}}
\newcommand{\bbC}{{\mathbb{C}}}

\newcommand{\calP}{{\mathcal P}}
\newcommand{\calC}{{\mathcal C}}

\newcommand{\lb}{\label}
\newcommand{\f}{\frac}

\newcommand{\bi}{\bibitem}

\newcommand{\beq}{\begin{equation}}
\newcommand{\eeq}{\end{equation}}
\newcommand{\ba}{\begin{align}}
\newcommand{\ea}{\end{align}}
\newcommand{\eps}{\varepsilon}
\newcommand{\ka}{\kappa}
\newcommand{\tht}{\theta}
\newcommand{\al}{\alpha}
\newcommand{\be}{\beta}
\newcommand{\ga}{\gamma}
\newcommand{\la}{\lambda}
\newcommand{\del}{\delta}
\newcommand{\til}{\tilde}

\newcounter{smalllist}

\DeclareMathOperator{\Real}{Re}

\allowdisplaybreaks
\numberwithin{equation}{section}
\newtheorem{theorem}{Theorem}[section]

\newtheorem*{p2.1}{Proposition 2.1}
\newtheorem{proposition}[theorem]{Proposition}
\newtheorem{lemma}[theorem]{Lemma}

\theoremstyle{definition}
\newtheorem{definition}[theorem]{Definition}

\newtheorem{conjecture}[theorem]{Conjecture}

\theoremstyle{remark}

\newcommand{\abs}[1]{\lvert#1\rvert}

\begin{document}

\title[Coefficients of OPUC and Higher-Order Szeg\H{o} Theorems]
{Coefficients of Orthogonal Polynomials on the Unit Circle and
Higher Order Szeg\H{o} Theorems}
%\author[L. Golinskii and A. Zlato\v{s}]{Leonid Golinskii$^{1}$ and Andrej Zlato\v{s}$^{2}$}

\author[L. Golinskii]{Leonid Golinskii$^*$}
\address{Mathematics Division, Institute for Low Temperature Physics and Engineering, 47 Lenin Avenue,
Kharkov 61103, Ukraine}
\email{golinskii@ilt.kharkov.ua}

\author[A. Zlato\v{s}]{Andrej Zlato\v{s}$^\dagger$}
\address{Department of Mathematics, University of Wisconsin, 480 Lincoln Drive, Madison, WI
53706, USA} \email{zlatos@math.wisc.edu}

\thanks{$^*$ The work of the first author was supported in part by INTAS Research Network NeCCA
03-51-6637.}
\thanks{$^\dagger$ Corresponding author}
%\thanks{$^1$ Mathematics Division, Institute for Low Temperature Physics and Engineering, 47 Lenin Avenue,
%Kharkov 61103, Ukraine; e-mail: {\tt golinskii@ilt.kharkov.ua}}

%\thanks{$^2$ Department of Mathematics, University of Wisconsin, 480 Lincoln Drive, Madison, WI 53706; e-mail:
%{\tt zlatos@math.wisc.edu}}

\begin{abstract}
Let $\mu$ be a non-trivial probability measure on the unit circle
$\partial\bbD$, $w$ the density of its absolutely continuous part,
$\al_n$ its Verblunsky coefficients, and $\Phi_n$ its monic
orthogonal polynomials. In this paper we compute the coefficients
of $\Phi_n$ in terms of the $\al_n$. If the function $\log w$ is
in $L^1(d\tht)$, we do the same for its Fourier coefficients. As
an application we prove that if $\al_n\in\ell^4$ and $Q(z) \equiv
\sum_{m=0}^N q_m z^m$ is a polynomial, then with $\bar Q(z) \equiv
\sum_{m=0}^N \bar q_m z^m$ and $S$ the left shift operator on
sequences we have
\[
\big| Q(e^{i\tht}) \big|^2\log w(\tht) \in L^1(d\tht) \,\,
\Leftrightarrow \,\,  \{\bar Q(S)\al\}_n\in\ell^2
\]
We also study relative ratio asymptotics of the reversed polynomials
$\Phi_{n+1}^*(\mu)/\Phi_n^*(\mu)-\Phi_{n+1}^*(\nu)/\Phi_n^*(\nu)$
and provide a necessary and sufficient condition in terms of the
Verblunsky coefficients of the measures $\mu$ and $\nu$ for this
difference to converge to zero uniformly on compact subsets of
$\bbD$.
\end{abstract}

\maketitle

%%%%%%%%%%%%%%%%%%%%%%%%%%%%%%%%%%%%%%%%%%%%%%%%%%%%%
\section{Introduction} \lb{S1}
%%%%%%%%%%%%%%%%%%%%%%%%%%%%%%%%%%%%%%%%%%%%%%%%%%%%%

In the present paper we study certain aspects of the theory of
orthogonal polynomials on the unit circle (OPUC). For background
information on the subject we refer the reader to the texts
\cite{Gbk,Sim1,Sim2,Szb}. Throughout, $d\mu$ will be a non-trivial
(i.e., with infinite support) probability measure on the unit
circle $\partial \bbD$ in $\bbC$, identified with the interval
$[0,2\pi)$ via the map $\tht\mapsto e^{i\tht}$. We will write
\[
d\mu(\tht) = w(\tht)\,\frac {d\tht}{2\pi} + d\mu_{\rm sing}(\tht)
\]
with $d\tht$ the Lebesgue measure on $[0,2\pi)$ and $d\mu_{\rm
sing}$ the singular part of $d\mu$.

One usually denotes by
\begin{equation} \lb{1.1}
\Phi_n(z) = \ka_{n,n}z^n + \ka_{n,n-1}z^{n-1} + \dots + \ka_{n,1}z
+ \ka_{n,0}
\end{equation}
the {\it monic} (i.e., $\ka_{n,n}=1$) {\it orthogonal polynomials}
for $\mu$ (with $n\ge 0$). It is standard to define the {\it
reversed polynomials} by
\begin{align*}
\Phi_n^*(z) & = \la_{n,n}z^n + \la_{n,n-1}z^{n-1} + \dots +
\la_{n,1}z + \la_{n,0}
\\ & \equiv \bar\ka_{n,0} z^n + \bar\ka_{n,1}z^{n-1} +
\dots + \bar\ka_{n,n-1}z + \bar\ka_{n,n}
\end{align*}
and let $\ka_{n,m}=\la_{n,m}=0$ whenever $m>n$. We have
$\Phi_0\equiv\Phi_0^*\equiv 1$ and for $n\ge 0$ the recurrence
relations
\begin{align}
\Phi_{n+1}(z) &= z\Phi_n(z) - \bar\al_n \Phi_n^*(z) \lb{1.2}
\\ \Phi_{n+1}^*(z) &= \Phi_n^*(z) - \al_n z \Phi_n(z) \lb{1.3}
\end{align}
with $\al_n\in\bbD$ the {\it Verblunsky coefficients} of $\mu$. A
fundamental result of Verblunsky \cite{V36} says that there is a
one-to-one correspondence between non-trivial probability measures
$\mu$ on $\partial\bbD$ and sequences $\{\al_n\}_{n\ge
0}\in\bbD^{\bbZ_0^+}$. If we set $\Phi_{n}\equiv 0$ and
$\Phi_{n}^*\equiv 1$ for $n\le -1$, and
%\begin{align*}
\[
\al_{-1} \equiv -1, \qquad \al_{n} \equiv 0 \quad (n\le -2)
\]
%\end{align*}
then \eqref{1.2}, \eqref{1.3} hold for all
$n\in\bbZ$. We accordingly let $\ka_{n,m}=0$ and
$\la_{n,m}=\del_{m,0}$ when $n<0$ and $m\ge 0$.

Probably the most famous OPUC result is {\it Szeg\H o's Theorem}. In
the form proved by Verblunsky \cite{V36} it says that
$\al_n\in\ell^2(\bbZ_0^+)$ if and only if $\log w(\tht)\in
L^1(d\tht)$.
%, in which case we say that $\mu$ is a {\it Szeg\H o measure}.
More precisely, the {\it sum rule}
\begin{equation} \lb{1.1a}
 \sum_{n=0}^\infty \log(1-\abs{\alpha_n}^2) = \int \log
(w(\theta)) \, \f{d\theta}{2\pi}
\end{equation}
holds. Note that both sides of \eqref{1.1a} are indeed non-positive
since $|\al_n|<1$ and by Jensen's inequality, $\int \log
(w(\theta))\, \f{d\theta}{2\pi}\leq \log (\int
w(\theta)\f{d\theta}{2\pi}) \leq \log (\mu(\partial\bbD))=0$, but
they can simultaneously be $-\infty$. Recently the area of sum
rules, for orthogonal polynomials as well as Schr\" odinger
operators, saw a rapid development starting with papers by
Deift-Killip \cite{DeiftK} and Killip-Simon \cite{KS}, which were
followed by many others (e.g.,
\cite{Den,Kup1,Kup2,LNS,NPVY,Sim,SZ2,Zla,Zla2}).

If $\al_n\in\ell^2$, one defines the {\it Szeg\H o function}
\[
D(z)\equiv \exp \bigg( \int \frac{e^{i\tht}+z}{e^{i\tht}-z} \log
w(\tht)\, \frac {d\tht}{4\pi} \bigg)
\]
which is analytic in $\bbD$. Szeg\H o's Theorem in its full extent
also shows that then
\begin{equation} \lb{1.1b}
\frac{\Phi_n^*}{\|\Phi_n^*\|_{L^2(d\mu)}} \to D^{-1}
\end{equation}
uniformly on compact subsets of $\bbD$. We have
\begin{equation} \lb{1.1c}
\|\Phi_n^*\|_{L^2(d\mu)} = \prod_{k=0}^{n-1} (1-|\al_k|)^{1/2} =
\prod_{k=0}^{n-1} \rho_k
\end{equation}
where $\rho_k\equiv \sqrt{1-|\al_k|^2}$ (see (1.5.13) in
\cite{Sim1}), and so if we define $d_m$ by
\begin{equation} \lb{1.3a}
D(z)^{-1} \equiv \bigg( \prod_{k\ge 0} \rho_k \bigg)^{-1}
(1+d_1z+d_2z^2+\dots)
\end{equation}
then
\begin{equation} \lb{1.4}
d_m = \lim_{n\to\infty} \la_{n,m}
\end{equation}

The first contribution of this paper is the following expression
of the coefficients $\ka_{n,m}$, $\la_{n,m}$, and $d_m$ in terms
of the $\al_k$. To the best of our knowledge (and to our
surprise), this result is new despite the long history and
classical nature of the subject!

\begin{theorem} \lb{T.1.1}
For $m\ge 1$,
\begin{equation} \lb{1.5}
\bar\ka_{n,n-m} = \la_{n,m} = \sum_{\substack{\sum_{1}^j a_l=m \\
j,a_l\ge 1}} \sum_{\substack{k_1< n \\ k_2< k_1-a_1 \\
\cdots \\ k_j<k_{j-1}-a_{j-1}}} \al_{k_1} \bar \al_{k_1-a_1} \dots
\al_{k_j} \bar \al_{k_j-a_j}
\end{equation}
If $\al_k\in\ell^2$, then also
\begin{equation} \lb{1.6}
d_m = \sum_{\substack{\sum_{1}^j a_l=m \\ j,a_l\ge 1}} \sum_{\substack{k_2< k_1-a_1 \\
\cdots \\ k_j<k_{j-1}-a_{j-1}}} \al_{k_1} \bar \al_{k_1-a_1} \dots
\al_{k_j} \bar \al_{k_j-a_j}
\end{equation}
\end{theorem}

{\it Remarks.} 1. In the above sums $[a_1,a_2,\dots,a_j]$ runs
through all $2^{m-1}$ ordered partitions of $m$, and $k_i\in\bbZ$.
%distinct variables and so, for instance, cases $[a_1,a_2]=[1,2]$ and
%$[a_1,a_2]=[2,1]$ are counted separately. Also, $j$ and $a_\ell$ are
%positive integers while $k_i\in\bbZ$.
\smallskip

2. Our choice of $\alpha_n$ for negative $n$ shows that the
condition ``$-1\le k_j-a_j$'' can be added under the second sum in
\eqref{1.5} (which is actually finite) and \eqref{1.6}. For
instance, if $m=n$, then the sum in \eqref{1.5} has a single
non-zero term with $j=1$, $a_1=n$, $k_1=n-1$, and so
$\ka_{n,0}=-\bar\al_{n-1}$. This can be seen from \eqref{1.2} and
$\Phi_n^*(0)=1$ as well.
\smallskip

3. Notice that for each partition $\{a_l\}$ with $\sum_{1}^j a_l=m$,
the second sum in \eqref{1.6} converges when $\al_k\in\ell^2$. This
is because then $\al_{k}\bar\al_{k-a}\in\ell^1$ for any fixed $a$,
and so
\[
|d_m|\le \sum_{\substack{\sum_{1}^j a_l=m \\ j,a_l\ge
1}}\prod_{l=1}^j \bigg( \sum_k |\al_{k}\bar\al_{k-a_l}| \bigg)
\]
%This will also be the case later.
\smallskip

Next, we describe an application of Theorem \ref{T.1.1} that
actually motivated our work. It involves the computation of  Taylor
coefficients of $\log D$. These are interesting primarily because
they coincide with Fourier coefficients of $\log w$. Indeed,
\begin{equation} \lb{1.9c}
\frac 12\, \frac{e^{i\tht}+z}{e^{i\tht}-z} = \frac 12 +
e^{-i\tht}z + e^{-2i\tht}z^2 + \dots
\end{equation}
and the definition of $D$ show that
\[
\log D(z) = \frac 12\, w_0 + w_1 z + w_2 z^2 +\dots
\]
where $w_m$ are defined by
\begin{equation} \lb{1.9b}
w_m \equiv \int e^{-im\tht} \log w(\tht) \frac{d\tht}{2\pi} = \bar
w_{-m}
%, \qquad \log w=\sum_{j=-\infty}^\infty w_j e^{ij\tht}
%\equiv w_0 + w_1 e^{i\tht} + \bar w_1 e^{-i\tht} + w_2 e^{2i\tht} + \bar w_2 e^{-2i\tht} + \dots
\end{equation}
We know from \eqref{1.1a} that
\begin{equation} \lb{1.9a}
w_0 = \sum_{k\ge 0} \log(1-\abs{\alpha_k}^2) = 2 \sum_{k\ge 0}
\log \rho_k
\end{equation}
and the methods from \cite{SZ2} can be used to compute the first few
of the other $w_m$. However,the corresponding computations become
very complicated with increasing $m$ (already at $m=4$ they are
close to intractable; \cite{SZ2} only deals with $m\le 2$). Our
method will provide $w_m$ for all $m$, although the resulting
formulae will obviously not be simple. That is why we postpone the
exact expressions to Theorem \ref{T.1.6} below and state the result
here in the following form that is sufficient for our first
application, Theorem \ref{T.3.3} (see also Lemma \ref{L.3.1} that
contains a similar formula for Taylor coefficients of $\log
\Phi_n^*$).

\begin{theorem} \lb{T.1.6x}
If $\al_k\in\ell^2$, then
\begin{equation} \lb{1.16}
w_m = \al_{m-1} - \sum_{k\ge 0} \al_{k+m} \bar\al_{k} + R_m(\mu)
\end{equation}
with
\begin{equation} \lb{1.17}
|R_m(\mu)| \le C_m \bigg( \sum_{k=0}^{m-1} |\al_k|^2 +
\sum_{k=m}^\infty |\al_k|^4 \bigg)
\end{equation}
\end{theorem}

%As we shall see, the bound \eqref{1.17} will play a crucial role
%in the proof of one of our main results, Theorem 1.4 below.

We note that \eqref{1.16} will be obtained from \eqref{1.6} by means
of expanding $\log(1+d_1z+d_2z^2+\dots)$ into its Taylor series.
This is a truly remarkable fact since the sum in \eqref{1.6} is
$m$-fold infinite and one might expect this method to only add
another degree of difficulty. Nevertheless, after appropriate
combinatorial manipulations it will turn out that the sum in
\eqref{1.16} (as well as the one in the exact form \eqref{1.15}) has
only a single infinite index!

The first application of the knowledge of $w_m$ we present in this
paper aims at the following conjecture of Simon \cite{Sim1} that is
a higher order generalization of \eqref{1.1a}. Here $S$ is the
left-shift operator on sequences
\begin{equation} \lb{1.9}
S(x_0,x_1,...)=(x_1,x_2,...)
\end{equation}

\begin{conjecture} \lb{C.1.2}
For distinct $\{\theta_m\}_{m=1}^l$ in $[0,2\pi)$ and $n_m$ positive
integers, define $N\equiv \sum_{m=1}^l n_m$, $n\equiv 1+\max_m n_m$,
and
%\begin{align*}
%Q(z) & \equiv \prod_{m=1}^\ell
%(z-e^{i\tht_m})^{\kappa_m}=\sum_{m=0}^N q_m z^m
%\\ \bar Q(z) & \equiv \prod_{m=1}^\ell (z-e^{-i\tht_m})^{\kappa_m}=\sum_{m=0}^N
%\bar q_m z^m
%\end{align*}
\begin{equation*}
Q(z)  \equiv \prod_{m=1}^l (z-e^{i\tht_m})^{n_m}=\sum_{m=0}^N q_m
z^m \qquad\text{and}\qquad \bar Q(z) \equiv \prod_{m=1}^l
(z-e^{-i\tht_m})^{n_m}=\sum_{m=0}^N \bar q_m z^m
\end{equation*}
so that
\[
\bar Q(S)=\sum_{m=0}^N \bar q_m S^m
\]
Then
\begin{equation} \lb{1.10}
%\begin{split}
\big| Q(e^{i\tht}) \big|^2\log w(\tht) \in L^1(d\tht) \,\,
\Leftrightarrow \,\, \{ \bar Q(S)\al\}_k\in\ell^2 \text{ and }
\al_k\in \ell^{2n}
%\end{split}
\end{equation}
\end{conjecture}

For $N=0$ this is just \eqref{1.1a}. For $N=1$ the conjecture was
proved by Simon (Theorem 2.8.1 in \cite{Sim1}) and for $N=2$ by
Simon and Zlato\v s \cite{SZ2}. It remains open for $N\ge 3$
although Denisov and Kupin \cite{DK}, mimicking the work of Nazarov,
Peherstorfer, Volberg, and Yuditskii \cite{NPVY} on Jacobi matrices,
showed that for each $Q$ there indeed is a condition in terms of
finiteness of a sum involving the $\al_k$ that is equivalent to the
LHS of \eqref{1.10}. Unfortunately, this sum is far from transparent
and its relation to the RHS of \eqref{1.10} is unclear.

Our contribution in this direction is the following {\it higher
order Szeg\H o theorem in $\ell^4$} which shows that Conjecture
\ref{C.1.2} holds if we a priori assume $\al_k\in\ell^4$.

\begin{theorem} \lb{T.3.3}
Assume that $\al_k\in\ell^4$, and for $q_0,q_1,\dots,q_N\in\bbC$
define
%\begin{align*}
%Q(z) & \equiv \sum_{m=0}^N q_m z^m
%\\ \bar Q(z) & \equiv \sum_{j=0}^N \bar q_m z^m
%\end{align*}
\begin{equation*}
Q(z) \equiv \sum_{m=0}^N q_m z^m \qquad \text{ and } \qquad \bar
Q(S)=\sum_{m=0}^N \bar q_m S^m
%\bar Q(z) \equiv \sum_{j=0}^N \bar q_m z^m
\end{equation*}
%so that
%\[
%\bar Q(S)=\sum_{m=0}^N \bar q_m S^m
%\]
Then
\begin{equation} \lb{3.10}
\big| Q(e^{i\tht}) \big|^2\log w(\tht) \in L^1(d\tht) \,\,
\Leftrightarrow \,\,  \{\bar Q(S)\al\}_k\in\ell^2
\end{equation}
\end{theorem}

{\it Remark.} Of course, the most interesting is the case from
Conjecture \ref{C.1.2} when all zeros of $Q$ are on the unit circle,
because the validity of the LHS of \eqref{3.10} only depends on
them.
\smallskip

Moreover, we provide in Theorem \ref{T.3.2} an exact formula for the
value of
\[
Z_Q(\mu)\equiv \int |Q(e^{i\tht})|^2\log w(\tht) \frac
{d\tht}{4\pi}
\]
in terms of the $\alpha_n$. Since $Z_Q$ is an entropy
\cite{KS,Sim1}, it is upper semi-continuous with respect to weak
convergence of measures, and so $Z_Q(\mu)\ge \limsup_n Z_Q(\mu_n)$
with $\mu_n$ the {\it Bernstein-Szeg\H o approximations} of $\mu$
having Verblunsky coefficients $\{\al_0,\dots,\al_n,0,0,\dots\}$. We
show in Proposition \ref{C.3.4} that, in fact, we always have
$Z_Q(\mu) = \lim_n Z_Q(\mu_n)$, including the case when both sides
are $-\infty$ (they cannot be $+\infty$ as each $Z_Q$ is bounded
above; see Section \ref{S3}).

Finally, we apply our method to the computation of the {\it
relative ratio asymptotics}
$\Phi_{n+1}^*(\mu)/\Phi_n^*(\mu)-\Phi_{n+1}^*(\nu)/\Phi_n^*(\nu)$
where $\Phi_n^*(\mu)$ and $\Phi_n^*(\nu)$ are the reversed
polynomials of measures $\mu$ and $\nu$, respectively.

\begin{theorem} \lb{T.4.1}
Let $\mu$ and $\nu$ be two non-trivial probability measures on
$\partial\bbD$. Let $\{\al_n(\mu)\}$ and $\{\al_n(\nu)\}$,
respectively, be their Verblunsky coefficients and let
$\Phi_n^*(\mu)$ and $\Phi_n^*(\nu)$, respectively, be their
reversed monic orthogonal polynomials. Then
\begin{equation} \lb{4.1}
\frac{\Phi_{n+1}^*(\mu)}{\Phi_n^*(\mu)} -
\frac{\Phi_{n+1}^*(\nu)}{\Phi_n^*(\nu)} \rightarrow 0
%\frac
%{\Phi_{n+1}^* / \Phi_{n}^*} {\Psi_{n+1}^* / \Psi_{n}^*} \to 1
\end{equation}
uniformly on compact subsets of $\bbD$ as $n\to\infty$ if and only
if for any $\ell\ge 1$
%there is $c_\ell$ such that
\begin{equation} \lb{4.2}
\lim_{n\to\infty} \big[ \al_n(\mu)\bar\al_{n-\ell}(\mu) -
\al_n(\nu)\bar\al_{n-\ell}(\nu) \big] = 0
\end{equation}
%Moreover, $H(z)\equiv 1$ if and only if all $c_\ell=0$.
\end{theorem}

As a corollary of Theorem \ref{T.4.1}, we provide a simple new proof
of the results of Khrushchev \cite{Khr} and Barrios and L\' opez
\cite{BL} on ratio asymptotics $\Phi_{n+1}^*/\Phi_n^*$ as
$n\to\infty$ of the reversed polynomials (Theorem \ref{T.4.2}), as
well as their generalization (Theorem \ref{T.4.3}).

The paper is organized as follows. Section \ref{S2} computes the
Taylor coefficients of $\Phi^*_n$ and $\log D$ in terms of the
Verblunsky coefficients and proves Theorems \ref{T.1.1} and
\ref{T.1.6x}. Section \ref{S3} introduces the {\it step-by-step sum
rules} (see \cite{KS,SZ1,SZ2}) and proves Theorem \ref{T.3.3}.
Section \ref{S4} proves Theorem \ref{T.4.1}.

%%%%%%%%%%%%%%%%%%%%%%%%%%%%%%%%%%%%%%%%%%%%%%%%%%%%%
\section{Coefficients of $\Phi^*_n(z)$ and $\log D(z)$ in Terms of Verblunsky Coefficients} \lb{S2}
%%%%%%%%%%%%%%%%%%%%%%%%%%%%%%%%%%%%%%%%%%%%%%%%%%%%%

We start with the proof of our first result, Theorem \ref{T.1.1}.

\begin{proof}[Proof of Theorem \ref{T.1.1}]
From \eqref{1.2} we have for $n\in \bbZ$ and $m\ge 0$,
\[
\ka_{n+1,m}=\ka_{n,m-1} - \bar\al_{n} \la_{n,m}
\]
with the convention $\ka_{n,-1}\equiv 0$. Substituting this
repeatedly into a similar equality obtained from \eqref{1.3}, we
get for $m\ge 1$,
\begin{align*}
\la_{n+1,m} & = \la_{n,m}-\al_n\ka_{n,m-1}
\\ & = \la_{n,m} + \al_n\bar\al_{n-1}\la_{n-1,m-1} - \al_n\ka_{n-1,m-2}
\\ & = ...
\\ & = \la_{n,m} + \sum_{a=1}^{m-1} \al_n\bar\al_{n-a} \la_{n-a,m-a}
+ \al_n\bar\al_{n-m}
\end{align*}
where in the last equality we have used $\ka_{n-m,-1}=0$ and
$\la_{n-m,0}=1$. If we now iterate this and note that
$\la_{m-l,m}=0$ for $m\ge 1$ and $l<0$, we have
\begin{align}
\la_{n+1,m} = & \sum_{k\le n} \sum_{a=1}^{m-1} \al_k\bar\al_{k-a}
\la_{k-a,m-a} + \sum_{k\le n} \al_k\bar\al_{k-m} \notag
\\ = & \sum_{k\le n} \sum_{a=1}^{m-1} \be_{k,a} \la_{k-a,m-a} + \sum_{k\le
n} \be_{k,m}  \lb{1.7}
\end{align}
with $\be_{k,a}\equiv \al_k\bar\al_{k-a}$. Of course, terms with
$k<0$ are zero.

We will prove \eqref{1.5} by induction on $n$. If $n\le 0$ and $m\ge
1$, then it obviously holds as in that case both sides are zero.
Assume therefore that \eqref{1.5} holds up to some $n$ and all $m\ge
1$. Then \eqref{1.7} gives
\begin{align*}
\la_{n+1,m} = & \sum_{k\le n} \sum_{a=1}^{m-1} \be_{k,a}
\sum_{\substack{\sum_{1}^j a_l=m-a \\ j,a_l\ge 1}} \sum_{\substack{k_1< k-a \\ k_2< k_1-a_1 \\
\cdots \\ k_j<k_{j-1}-a_{j-1}}} \be_{k_1,a_1} \dots \be_{k_j,a_j} +
\sum_{k\le n} \be_{k,m}
\\ = &  \sum_{\substack{\sum_{0}^j a_l=m \\ j,a_l\ge 1}} \sum_{\substack{k_0< n+1 \\ k_1< k_0-a_0 \\
\cdots \\ k_j<k_{j-1}-a_{j-1}}} \be_{k_0,a_0} \dots \be_{k_j,a_j} +
\sum_{k_0< n+1} \be_{k_0,m}
\\ = &  \sum_{\substack{\sum_{0}^j a_l=m \\ j\ge
0\\ a_l\ge 1}} \sum_{\substack{k_0< n+1 \\ k_1< k_0-a_0 \\
\cdots \\ k_j<k_{j-1}-a_{j-1}}} \be_{k_0,a_0} \dots \be_{k_j,a_j}
\end{align*}
with $k_0\equiv k$ and $a_0\equiv a$. But this is \eqref{1.5} for
$n+1$ in place of $n$. Thus \eqref{1.5} is proved, and \eqref{1.6}
follows from \eqref{1.4}.
\end{proof}

Our next aim is to compute the Taylor coefficients of $\log D$. We
will again assume $\al_k\in\ell^2$ so that $D$ is well defined. By
\eqref{1.3a} we have for $z$ close to 0,
\[
\log \bigg( \prod_{k\ge 0} \rho_k \bigg) - \log D(z) = \sum_{j\ge 1}
\frac {(-1)^{j-1}} j (d_1 z + d_2 z^2 +\dots)^j
\]
and so $w_m$ is the negative of the $m^{\rm th}$ Taylor coefficient
of the RHS when $m\ge 1$. That is,
\begin{align}
w_m & = \sum_{\substack{\sum_{1}^j b_\ell=m \\ j,b_\ell\ge 1}} \frac
{(-1)^{j}} j \prod_{\ell=1}^j d_{b_\ell}  \notag
\\ & = \sum_{\substack{\sum_{1}^j b_\ell=m \\ j,b_\ell\ge 1}} \frac
{(-1)^{j}} j \prod_{\ell=1}^j \bigg( \sum_{\substack{\sum_{1}^{p}
a_l=b_\ell \\ p,a_l\ge 1}} \sum_{\substack{k_2< k_1-a_1 \\
\cdots \\ k_p<k_{p-1}-a_{p-1}}} \be_{k_1,a_1} \dots \be_{k_p,a_p}
\bigg)  \lb{1.11}
\\ & = \sum_{\{ (k_1,a_1),\dots ,(k_i,a_i) \}\in M_m} \be_{k_1,a_1}
\dots \be_{k_i,a_i}  \sum_{j=1}^i \frac {(-1)^{j}} j N_j
\big( \{ (k_1,a_1),\dots ,(k_i,a_i) \} \big)  \lb{1.12}
\end{align}
%\bigg[ \prod_{l=1}^i \be_{k_l,a_l} \bigg]
with $M_m$ and $N_j$ defined below.

Before stating the definitions, let us first describe how
\eqref{1.12} was obtained from \eqref{1.11}. We multiply out the
brackets in \eqref{1.11} to get a sum of products $\be_{k_1,a_1}
\dots \be_{k_i,a_i}$ (with coefficients), and then collect terms
with identical products (only differing by a permutation). The
coefficient at each product $\be_{k_1,a_1} \dots \be_{k_i,a_i}$
obtained in this way will then equal the last sum in \eqref{1.12}.
For example, the product $\be_{3,1}\be_{1,1}^2$ appears in
\eqref{1.11} for $m=3$ as $\tfrac{(-1)^3}3
(\be_{3,1})(\be_{1,1})(\be_{1,1})$, $\tfrac{(-1)^3}3
(\be_{1,1})(\be_{3,1})(\be_{1,1})$, $\tfrac{(-1)^3}3
(\be_{1,1})(\be_{1,1})(\be_{3,1})$, $\tfrac{(-1)^2}2
(\be_{3,1}\be_{1,1})(\be_{1,1})$, and $\tfrac{(-1)^2}2
(\be_{1,1})(\be_{3,1}\be_{1,1})$. The first three come from $j=3$
and $b_1=b_2=b_3=1$ in \eqref{1.11}, the fourth from $j=2$, $b_1=2$,
$b_2=1$, and the fifth from $j=2$, $b_1=1$, $b_2=2$. Therefore the
coefficient at $\be_{3,1}\be_{1,1}^2$ in \eqref{1.12} has to be
$\tfrac{(-1)^3}3 3+\tfrac{(-1)^2}2 2=0$.

It is obvious that the products that appear in \eqref{1.12} must
satisfy $i,a_l\ge 1$ and $\sum_1^i a_l = m$, because the sum of the
$a_l$'s in any term of the $\ell^{\rm th}$ bracket of \eqref{1.11}
equals $b_\ell$. The set $M_m$ will therefore reflect this
condition. The question now is, given any collection (i.e., set with
repetitions; see below) of couples $P=\{ (k_1,a_1),\dots ,(k_i,a_i)
\}\in M_m$, in how many ways can the corresponding product
$\be_{k_1,a_1} \dots \be_{k_i,a_i}$ be obtained by multiplying out
$j\ge 1$ brackets in \eqref{1.11}. If this number is denoted
$N_j'(P)$, then the correct coefficient at $\be_{k_1,a_1} \dots
\be_{k_i,a_i}$ in \eqref{1.12} is $\sum_{j=1}^i \frac {(-1)^{j}} j
N_j'(P)$. Hence to obtain \eqref{1.11}=\eqref{1.12}, we are left
with showing that $N_j(P)$, defined below, equals $N_j'(P)$.

We will call a {\it collection} an unordered list of elements,
some of which can be identical (i.e., a collection is a set that
can contain multiple identical elements, a hat with multicolored
balls). Such identical elements are considered indistinguishable.
A $j$-{\it tuple} will be an ordered list of $j$ elements.
Collections will be denoted by $\{\dots\}$, $j$-tuples by
$[\dots]$. Below we will consider collections and $j$-tuples whose
elements are couples $(k,a)$ with $k\in\bbZ$, $a\in\bbN$. For
instance, $\{(3,1),(1,1),(1,1)\}$ is a collection
($\{(1,1),(1,1),(3,1)\}$ is the same one) and
$[(3,1),(1,1),(1,1)]$, $[(1,1),(3,1),(1,1)]$,
$[(1,1),(1,1),(3,1)]$ are three distinct triples. Finally, the
union of collections is the collection obtained by joining their
lists of elements, for instance, $\{(3,1),(1,1)\}\cup\{(1,1)\} =
\{(3,1),(1,1),(1,1)\}$.

\begin{definition} \lb{D.1.3}
Let $M_m$ be the set of all distinct collections
$P=\{(k_1,a_1),\dots ,(k_i,a_i)\}$ with $i\ge 1$, $k_l\in\bbZ$, and
$a_l\ge 1$ such that $\sum_1^i a_l = m$. We let
\begin{equation} \lb{1.12a}
\be(P) \equiv \be_{k_1,a_1}\dots \be_{k_i,a_i} =
\al_{k_1}\bar\al_{k_1-a_1}\dots \al_{k_i}\bar\al_{k_i-a_i}
\end{equation}

We say that a collection $P=\{(k_1,a_1),\dots ,(k_i,a_i)\}$ is {\it
linear} if $k_u<k_v-a_v$ or $k_v<k_u-a_u$ whenever $u\neq v$. In
particular, $k_u\neq k_v$ when $u\neq v$, which means that a linear
collection $P$ cannot contain two identical couples, and thus it is
just a set.
%A collection $P=\{(k_1,a_1)\}$ is always linear.

If $P\in M_m$ and $j\ge 1$, then $N_j \big( P\big)$ is the number of
distinct $j$-tuples $\calP=[P_1,\dots,P_j]$ such that each $P_\ell$
is a non-empty linear collection and $\bigcup_{\ell=1}^j P_\ell =
P$. We will call each such $\calP$ an {\it admissible division} of
$P$.
\end{definition}

%In particular, linearity of $P_\ell$ means that two identical
%couples $(k_l,a_l)$ from $P$ cannot be in the same $P_\ell$ (and
%hence $P_\ell$ is just a set). Also notice that since identical
%couples are indistinguishable, divisions $\calP$ differing only by a
%permutation of identical couples are not considered distinct.
For instance, if $P=\{ (3,1),(1,1),(1,1)\}$ (corresponding to
$\be_{3,1}\be_{1,1}^2$ above), then the admissible divisions are
$[\{(3,1)\},\{(1,1)\},\{(1,1)\}]$,
$[\{(1,1)\},\{(3,1)\},\{(1,1)\}]$, $[\{(1,1)\},\{(1,1)\},\{(3,1)\}]$
(with $j=3$) and $[\{(3,1),(1,1)\},\{(1,1)\}]$,
$[\{(1,1)\},\{(3,1),(1,1)\}]$ (with $j=2$). Hence in this case
$N_3(P)=3$, $N_2(P)=2$, $N_1(P)=0$ and the last sum in \eqref{1.12}
is indeed $\tfrac{(-1)^3}3 3+\tfrac{(-1)^2}2 2=0$.

To finish the proof of \eqref{1.11}=\eqref{1.12} we need to show
that $N_j(P)=N_j'(P)$ for any $P\in M_m$ (as we did for $P=\{
(3,1),(1,1),(1,1)\}$), where $N_j'(P)$ is the number of times the
product $\be(P)=\be_{k_1,a_1} \dots \be_{k_i,a_i}$ is obtained by
multiplying out $j$ brackets in \eqref{1.11}. The desired equality
follows from realizing that the collection $P_\ell$ in the
definition ($\ell=1,\dots,j$) corresponds to the ``subproduct'' of
$\be(P)$ coming from the $\ell^{\rm th}$ bracket in \eqref{1.11}
(which is why the $\calP$'s must be ordered, as well as why
$P_\ell$ must be linear). With this identification in mind, it is
easy to see that each admissible division $[P_1,\dots,P_j]$ of $P$
corresponds to precisely one way of obtaining $\be(P)$ in
\eqref{1.12} from \eqref{1.11} by multiplying $j$ subproducts
(from $j$ brackets) corresponding to $P_1,\dots,P_j$ (with
$b_\ell$ being the sum of the $a_l$ for which $(k_l,a_l)\in
P_\ell$), and vice versa.

%Note that order {\it within} each subproduct as well as {\it within}
%each $P_\ell$ is already determined by the common condition
%$k_u<k_v-a_v$; that is why each subset $P_\ell$ is unordered. This
%ensures that for each $P=\{(k_1,a_1),\dots ,(k_i,a_i)\}\in M_m$ the
%sum $\sum_{j=1}^i \tfrac {(-1)^{j}} j N_j \big( P \big)$ in
%\eqref{1.12} equals the coefficient at the product $\be_{k_1,a_1}
%\dots \be_{k_i,a_i}$ in \eqref{1.11} with the brackets multiplied
%out.

Hence we have obtained an explicit expression for $w_m$. We will now
simplify it considerably by showing that coefficients at many
$\be(P)$ in \eqref{1.12} are actually zero, as was the case for
$\be_{3,1}\be_{1,1}^2$ (see Lemma \ref{L.1.5} below).

We say that $K\in\bbZ$ is a {\it cut} of $P=\{ (k_1,a_1),\dots
,(k_i,a_i) \}\in M_m$ if $k_l\neq K$ for all $l$, if $\min_l
\{k_{l}\}<K<\max_l \{k_{l}\}$, and if for any $u,v$ with
$k_{u}<K<k_{v}$ we have $k_u<k_v-a_v$. For instance, $P=\{
(3,1),(1,1),(1,1)\}$ has one cut $K=2$. Our interest here will be
mainly in ``cuttless'' collections as is demonstrated by the
following two lemmas.

\begin{lemma} \lb{L.1.4}
If $P\in M_m$ has no cut, then
\begin{equation} \lb{1.12b}
\max_l \{k_l\} - \min_l \{k_l-a_l\} \le m
\end{equation}
\end{lemma}

\begin{proof}
Consider the union of intervals $I\equiv\bigcup_l
[k_l-a_l,k_l]\subset\bbR$, with $|I|\le \sum_{l} a_l =m$. If $P$ has
no cut, then $I$ is an interval (and vice versa) because otherwise
the minimum of any component, except for the bottom one, were a cut.
But then we obviously have
\[
I=[\min_l \{k_l-a_l\}, \max_l \{k_l\}]
\]
proving \eqref{1.12b}.
\end{proof}

Let $|P|$ be the number of elements of a collection $P$, counting
identical elements as many times as they are included in $P$. For
instance, $|\{(3,1),(1,1),(1,1)\}|=3$.

\begin{lemma} \lb{L.1.5}
If $P\in M_m$ has a cut, then
\begin{equation} \lb{1.13}
\sum_{j=1}^{|P|} \frac {(-1)^{j}} j N_j \big( P \big) = 0
\end{equation}
\end{lemma}

\begin{proof}
Fix $P=\{ (k_1,a_1),\dots ,(k_i,a_i) \}\in M_m$ that has a cut $K$.
First notice that if $\Pi$ is the set of all admissible divisions
$\calP=[P_1,\dots,P_{j_\calP}]$ of $P$, then \eqref{1.13} is equivalent to
\begin{equation} \lb{1.13a}
\sum_{\calP\in\Pi} \frac {(-1)^{j_\calP}} {j_\calP} = 0
\end{equation}

For each $\calP$ let $\calC(\calP)$ be the collection (not a union!) of up to $2j_\calP$
non-empty sets that we obtain by splitting each $P_\ell$ at $K$.
That is, we define
\[
P_\ell^\pm \equiv \{(k_l,a_l)\in P_\ell \,\big|\, \pm (k_l-K)>0\}
\]
%(recall that $P_\ell$ is just a set, and hence so are $P_\ell^\pm$),
so that $P_\ell=P_\ell^+\cup P_\ell^-$, and then let $\calC(\calP)$
be the collection of those $P_\ell^\pm$ that are not empty.
%whenever $P_\ell$ contains couples $(k_u,a_u), (k_v,a_v)$ with
%$k_u<K<k_v$, we split it into two subsets of those elements with
%$k_l$ smaller/larger than $K$, and then take the collection of (up
%to $2j_\calP$) so obtained sets. Notice that these are indeed sets
Notice that the $P_\ell^\pm$ are indeed sets and they are linear ---
both because the same is true for $P_\ell$.

%Consider the collection of all divisions $\calP$ of $P$ into subsets
%$P_1,\dots,P_{j_\calP}$ described above (for all $j_\calP$) and
%divide it into equivalence classes, in each divisions with the same
%$\calC(\calP)$. We will show that for each such class, with
%if $\calC_0=\calC(\calP)$ for some $\calP$, then
Hence $\calC$ defines an equivalence relation on $\Pi$ by
$\calP\sim\calP'$ iff $\calC(\calP)=\calC(\calP'$). We will show
that the part of the sum in \eqref{1.13a} corresponding to any
equivalence class is zero. That is, we will prove
\begin{equation} \lb{1.14}
\sum_{\calC(\calP)=\calC_0} \frac {(-1)^{j_\calP}}{j_\calP} = 0
\end{equation}
for any $\calC_0$ such that $\calC(\calP)=\calC_0$ for some $\calP\in\Pi$.

Let us fix any such $\calC_0$. Then $\calC_0$ is a collection of
non-empty linear sets $Q_1,\dots,Q_q$ and $R_1,\dots,R_r$ whose
union (as a union of collections) is $P$, such that if
$(k_l,a_l)\in Q_u$, then $k_l<K$, and if $(k_l,a_l)\in R_u$, then
$k_l>K$. That is, the $Q_u$ are the non-empty $P_\ell^-$ and the
$R_v$ are the non-empty $P_\ell^+$. Let $q\le r$, since the case
$q\ge r$ is identical.

Assume first that these sets are all distinct. Then for every $0\le
s\le q$ there are ${q \choose s}{r \choose s} s!(q+r-s)!$ admissible
divisions $\calP$ of $P$ with $\calC(\calP)=\calC_0$ and
$j_\calP=q+r-s$. These are created by choosing $s$ sets from
$Q_1,\dots,Q_q$ and $s$ from $R_1,\dots,R_r$, taking all $s!$
pairings of the selected $Q$'s with the selected $R$'s, and then all
$(q+r-s)!$ orderings of thus created $q+r-s$ sets (unions of the
paired couples $Q_u\cup R_v$ together with the
unpaired $Q$'s and $R$'s) --- the $P_\ell$'s. Since all the original
sets were distinct, this construction gives no repetitions. Notice
also that any $P_\ell=Q_u\cup R_v$ is linear because so are $Q_u$
and $R_v$ and $K$ is a cut for $P$. This shows that the LHS of
\eqref{1.14} equals
\[
\sum_{s=0}^q \frac{(-1)^{q+r-s}}{q+r-s} {q \choose s}{r \choose s}
s!(q+r-s)! = (-1)^{q+r}(q+r-1)! \sum_{s=0}^q (-1)^{s} \frac {{q
\choose s}{r \choose s}} {{q+r-1 \choose s}} = 0
\]
The last equality follows from Lerch's identity \cite{Ler} (also
in \cite [p. 61] {Kau})
\[
\sum_{s=0}^q (-1)^{s} \frac {{q \choose s}{r \choose s}} {{p \choose
s}} = \frac{{p-r \choose q}}{{p\choose q}}
\]
which holds whenever $p\ge q$.

If now some $Q$'s and/or some $R$'s are identical, then in the above
sum every $\calP$ with $\calC(\calP)=\calC_0$ is counted the same
number of times $T$, which equals the product of the factorials of
the numbers of identical sets. This is because there are $T$
permutations of the $Q$'s and $R$'s that fix the classes of
identical sets, and hence when we perform the above algorithm to
obtain all admissible $\calP$'s with $\calC(\calP)=\calC_0$, each
such $\calP$ will be obtained $T$ times. Therefore the LHS of
\eqref{1.14} equals
\[
\frac 1T \sum_{s=0}^q \frac{(-1)^{q+r-s}}{q+r-s} {q \choose s}{r
\choose s} s!(q+r-s)! = 0
\]
This proves \eqref{1.14}, and \eqref{1.13a} follows by summing over
all $\calC_0$.
\end{proof}

Hence the only terms that matter in \eqref{1.12} are those with no
cuts (which is the main point of this section). Moreover, it is
obvious that $\be(P)=0$ when some $k_l-a_l\le -2$. Therefore we
define $\omega(P)\equiv \max_l\{k_l \,|\, (k_l,a_l)\in P\}$,
$\del(P)\equiv \min_l\{k_l-a_l \,|\, (k_l,a_l)\in P\}$,
\begin{equation} \lb{1.14b}
N(P)\equiv \sum_{j=1}^{|P|} \frac {(-1)^{j}} j N_j (P)
\end{equation}
and for $0\le n\le\infty$
\begin{equation} \lb{1.14a}
M^n_m\equiv \{ P\in M_m \,\big|\, \text{$P$ has no cuts and
$0\le\omega(P)\le n$} \}
\end{equation}
If now $P\in M_m\smallsetminus M^\infty_m$, then either $P$ has a
cut and so $N(P)=0$, or $\omega(P)\le -1$ and then $\be(P)=0$
because $\delta(P)\le -2$. This means that the sum in \eqref{1.12}
only needs to be taken over $M_m^\infty$. Before formally stating
this fact, we remark that
\begin{equation} \lb{1.15a}
M^n_m=\{P+k \,\big|\, P\in M^0_m \text{ and } 0\le k\le n\}
\end{equation}
where $P+k\equiv \{(k_l+k,a_l)\,|\, (k_l,a_l)\in P\}$. Also notice
that $N(P+k)=N(P)$ by definition, $P-\omega(P)\in M^0_m$ for any
$P\in M_m^\infty$, and $M^0_m$ is a finite set by Lemma \ref{L.1.4}.
In this light the following result is an immediate consequence of
\eqref{1.12} and Lemma \ref{L.1.5}.

\begin{theorem} \lb{T.1.6}
If $\al_k\in\ell^2$, then for $m\ge 1$
\begin{equation} \lb{1.15}
w_m  = \sum_{P\in M^\infty_m} N(P) \be(P) = \sum_{P\in M^0_m} N(P)
\sum_{k=0}^\infty \be(P+k)
\end{equation}
%In particular, Theorem \ref{T.1.6x} holds.
\end{theorem}

%{\it Remark.} Here $O(\dots)$ obviously depends on $m$ and
%$\|\al\|_4^4$ includes $\al_{-1}=-1$. The proof shows that if we
%take $\|\al\|_p^p = \sum_{k=0}^\infty |\al_k|^p$, then
%$O(\|\al\|_4^4)$ has to be replaced by $O(\|\al\|_4^4 +
%\sum_{k=0}^{m-1} |\al_k|^2)$.

{\it Remark.} The second form of $w_m$ in \eqref{1.15} shows that
for $m\neq 0$, the $m^{\rm th}$ Fourier coefficient of $\log
w(\tht)$ (and so the $m^{\rm th}$ Taylor coefficient of $\log
D(z)$) can be expressed as a sum over a single infinite index of
products involving only ``nearby'' $\alpha_k$'s.
\smallskip

%2. As we shall see, the bound \eqref{1.17} will play a crucial
%role in the proof Theorem \ref{T.3.3}.

Now we are ready to prove Theorem \ref{T.1.6x}.

\begin{proof}[Proof of Theorem \ref{T.1.6x}]
%Equation \eqref{1.15} follows from \eqref{1.12} and Lemma
%\ref{L.1.5}. If we let
If
\[
M'_m\equiv \{ P\in M^\infty_m \,\big|\, |P|=1 \}
\]
then \eqref{1.15} can be written as
\begin{equation} \lb{1.18}
w_m= \bigg(\sum_{P\in M_m'} + \sum_{P\in M^\infty_m\smallsetminus
M_m'} \bigg) N(P) \be(P)
%\sum_{j=1}^{|P|} \frac {(-1)^{j}} j N_j(P)
\end{equation}
Note that the sum in \eqref{1.16} together with
$\al_{m-1}=-\al_{m-1}\al_{-1}$ is just the first sum in
\eqref{1.18}, and so $R_m(\mu)$ is the second sum in \eqref{1.18}.
It remains to prove \eqref{1.17}.

Since each $N_j(P)$ is bounded by a constant only depending on $m$,
\[
|R_m(\mu)|\le C_m \sum_{P\in M^\infty_m\smallsetminus M_m'}
|\be(P)|
\]
For any $P\in M^\infty_m\smallsetminus M_m'$ let  $i\equiv |P|\ge
2$. If $\del(P)\le -2$, then $\be(P)=0$. If $\del(P)=-1$, then by
Lemma \ref{L.1.4}
\[
|\be(P)|\le \prod_{l=1}^{i} |\al_{k_l}| \le \sum_{j=0}^{m-1}
|\al_j|^{i} \le \sum_{j=0}^{m-1} |\al_j|^2
\]
And if $\del(P)\ge 0$, then
\[
|\be(P)|\le \sum_{l=1}^{i} \big(
|\al_{k_l}|^{2i}+|\bar\al_{k_l-a_l}|^{2i} \big) \le m
\sum_{j=\del(P)}^{\del(P)+m} |\al_j|^{2i} \le m
\sum_{j=\del(P)}^{\del(P)+m} |\al_j|^{4}
\]
Since each $P\in M^\infty_m$ has no cuts, the number of $P\in
M^\infty_m$ with any given $\del(P)$ is a finite constant only
depending on $m$. Hence \eqref{1.17} follows and the proof is
complete.
\end{proof}

We write here explicitly the first three $w$'s from \eqref{1.15}.
Recall that $\al_{-1}=-1$, $\al_{-2}=\al_{-3}=\dots=0$, and
$\rho_k=\sqrt{1-|\al_k|^2}\,$.
\begin{align*}
w_1 =& -\sum_k \al_k\bar\al_{k-1}
\\ w_2 =& -\sum_k \al_k\bar\al_{k-2}\rho_{k-1}^2 + \tfrac 12
\sum_k \al_k^2\bar\al_{k-1}^2
\\ w_3 =& -\sum_k \al_k\bar\al_{k-3}\rho_{k-1}^2\rho_{k-2}^2 +
\sum_k \al_k^2\bar\al_{k-1}\bar\al_{k-2}\rho_{k-1}^2  + \sum_k
\al_k\al_{k-1}\bar\al_{k-2}^2\rho_{k-1}^2 - \tfrac 13 \sum_k
\al_k^3\bar\al_{k-1}^3
\end{align*}

Finally, we note that all Taylor coefficients of $\log (D(z)/D(0))$
verify the claim of the remark after Theorem \ref{T.1.6}. It turns
out that this is essentially the only such function of the form
$F(D(z)/D(0))$.

\begin{proposition} \lb{P.1.7}
Assume that $F$ is analytic on a neighborhood of $1$ and each
Taylor coefficient $h_{m}$ of $H(z)=F(D(z)/D(0))$ is, as a
function of $\{\al_k\}\in\ell^2$, a sum of products of the
$\alpha_k$'s such that if $\al_k$ and $\al_l$ both appear in the
same product, then $|k-l|\le c_{F,m}$ for some $c_{F,m}<\infty$.
It follows that $H(z)=a+b\log (D(z)/D(0))$ for some $a,b\in\bbC$.
\end{proposition}

\begin{proof}
Define $G(z)=F(e^z)$ so that $G$ is analytic on a neighborhood of
$0$ (with Taylor coefficients $g_m$) and $H(z)=G(\log
(D(z)/D(0)))$. The fact that $\log (D(z)/D(0)) = \sum_{m\ge 1}w_m
z^m$ satisfies the proposition shows that when $g_m$ is the first
non-zero coefficient with $m\ge 2$, then $h_m$ does not satisfy
the required condition because $h_m=g_1w_m+g_mw_1^m$. Therefore
$G(z)=a+bz$ for some $a,b$.
\end{proof}

%%%%%%%%%%%%%%%%%%%%%%%%%%%%%%%%%%%%%%%%%%%%%%%%%%%%%
\section{A Higher Order Szeg\H o Theorem} \lb{S3}
%%%%%%%%%%%%%%%%%%%%%%%%%%%%%%%%%%%%%%%%%%%%%%%%%%%%%

In this section we will prove Theorem \ref{T.3.3}. We will do this
by first deriving sum rules \` a la Denisov-Kupin \cite{DK} that
provide us a necessary and sufficient condition for the left hand
side \eqref{3.10} to hold. The difference between our Theorem
\ref{T.3.2} below and \cite{DK} is that in \cite{DK} this condition
is expressed in terms of traces of powers of the CMV matrix (see,
e.g., \cite{Sim1}), which is less explicit than the form we obtain
here (although, obviously, the two conditions have to be
equivalent). This, together with Theorem \ref{T.1.6}, will suffice
to yield Theorem \ref{T.3.3}.

We start by introducing some notation. The Carath\'eodory and
Schur functions, $F:\bbD\to i\bbC^-$ and $f:\bbD\to\bbD$, for
$d\mu$ are defined by
\[%\begin{align*}
F(z) \equiv \int \f{e^{i\theta}+z}{e^{i\theta}-z}\, d\mu(\theta)
\equiv \f{1+zf(z)}{1-zf(z)}
\] %\lb{1.3}
%\end{align*}
It is a result of Geronimus \cite{Ger44} that the Verblunsky
coefficients of $\mu$ coincide with the Schur parameters of $f$
defined inductively by the Schur algorithm
\begin{equation} \lb{3.0x}
f(z) = \f{\alpha_0 + zf_1 (z)}{1+z\bar\alpha_0 f_1(z)}
\end{equation}
Here \eqref{3.0x} defines $\alpha_0\in \bbD$ and $f_1:\bbD\to\bbD$,
and iteration then yields $\alpha_1, \alpha_2, \dots$ and $f_2, f_3,
\dots$. Note that $f(0)=\al_0$ and, by induction, $m^{\rm th}$
Taylor coefficient of $f$ only depends on $\al_0,\dots,\al_m$.

In the following we will write $\Phi_n^*(\mu,z)$ and $D(\mu,z)$
for the reversed polynomials and the Szeg\H o function.
Accordingly, we will write $w_m(\mu)$ for the Taylor coefficients
of $\log D(\mu,z)$, and we will also let
\[
-\log \Phi_n^*(\mu,z) \equiv \sum_{m\ge 1} w_{n,m}(\mu)z^m
\]

We will now fix a measure $\mu$ and denote its Verblunsky
coefficients $\al_k$. For the sake of transparency, we will include
$\al_{-1}=-1$ at the beginning of the sequence of the coefficients,
so that these will be $\{-1,\al_0,\al_1,\dots\}$. We let $\mu_n$ be
the $n$-th Bernstein-Szeg\H o approximation of $\mu$, with
Verblunsky coefficients
$\{-1,\al_0,\al_1,\dots,\al_{n},0,0,\dots\}$, and
$\mu^{(n)}=w^{(n)}(\tht)\tfrac{d\tht}{2\pi} +d\mu^{(n)}_{\rm sing}$
the measure with Verblunsky coefficients
$\{-1,\al_n,\al_{n+1},\dots\}$.
%obtained by deleting the zeroth through $(n-1)$-st coefficient of
%$\mu$.
%The corresponding Szeg\H o functions will be $D_j$ and
%$D^{(j)}$ (the latter is defined only when $D$ is). The polynomials
%$\Phi_n^*(\mu^{(j)})$ will be denoted $\Phi_n^{(j)*}$.
%We also let $w_{j;m}$, $w_m^{(j)}$, and $\del w_m$ be the $m$-th
%Taylor coefficients of $\log D_j(z)$, $\log D^{(j)}(z)$, and $\log
%\del D(z)$. The $m$-th Taylor coefficients of $\log \Phi^*_n(z)$
%and $\log \Phi^{(j)*}_n(z)$ will be $u_{n,m}$ and $u_{n,m}^{(j)}$.

In Section~2.9 of \cite{Sim1}, Simon defines the {\it relative
Szeg\H{o} function}
\begin{equation} \lb{3.0}
(\delta D)(\mu,z) \equiv \f{1-\bar\alpha_0 f(z)}{\rho_0} \,
\f{1-zf_1(z)}{1-zf(z)}
\end{equation}
with $f,f_1$ from \eqref{3.0x}. Its advantage is that, unlike $D$,
it is defined for any $\mu$. If $w(\tht)$ is positive almost
everywhere, then so is $w^{(1)}(\tht)$, and
\begin{equation} \lb{3.0z}
\log\frac{w(\tht)}{w^{(1)}(\tht)}\in L^p[0,2\pi), \quad p<\infty
\end{equation}
with
\begin{equation} \lb{3.0a}
(\delta D)(\mu,z) = \exp \bigg( \int \frac{e^{i\tht}+z}{e^{i\tht}-z}
\log \bigg( \frac{w(\tht)}{w^{(1)}(\tht)} \bigg) \, \frac
{d\tht}{4\pi} \bigg)
\end{equation}
(see Theorem 2.9.3 in \cite{Sim1}). Obviously, in the case
$\al_k\in\ell^2$ we have
\begin{equation} \lb{3.0b}
(\delta D)(\mu,z) = \frac{D(\mu,z)}{D(\mu^{(1)},z)}
\end{equation}
which explains the name. We define the Fourier coefficients of
$\log(w(\tht)/w^{(1)}(\tht))$ to be $\del w_m(\mu)$ so that from
\eqref{3.0a} and \eqref{1.9c} we obtain
\begin{equation} \lb{3.0d}
\log (\del D)(\mu,z) = \frac 12 \del w_0(\mu) + \del w_1(\mu)z +
\del w_2(\mu) z^2 +\dots
\end{equation}
In particular,
\begin{equation} \lb{3.0e}
\del w_0(\mu) = 2 \log\rho_0
\end{equation}
by \eqref{3.0} and $f(0)=\al_0$, and $\del w_{-m}(\mu) =
\overline{\del w_m(\mu)}$. In the case $\al_k\in\ell^2$ we have by
\eqref{3.0b}
\begin{equation} \lb{3.0c}
\del w_m(\mu) = w_m(\mu)-w_m(\mu^{(1)})
\end{equation}

Finally, we define $N(P)$ by \eqref{1.14b}, $\be(P)$ by
\eqref{1.12a}, and let $\be^{(n)}(P)$ be defined as $\be(P)$, but
with $\al_{-1},\al_0,\dots,\al_{n-2}$ replaced by zeros and
$\al_{n-1}$ replaced by $-1$. That is, $\be^{(n)}(P)$ equals
$\be(P-n)$ for the measure $\mu^{(n)}$. For instance,
$\be^{(2)}(\{(3,1),(1,1)\})=\al_3\bar\al_2(-1)\bar 0=0,$ which is
$\be(\{(1,1),(-1,1)\})$ for the measure $\mu^{(2)}$ with Verblunsky
coefficients $\{-1,\al_2,\al_3,\dots\}$. In particular, \eqref{1.15}
for the measure $\mu^{(n)}$ and $N(P-n)=N(P)$ imply
\begin{equation} \lb{3.1a}
w_m(\mu^{(n)}) =  \sum_{P\in M^\infty_m} N(P)\be^{(n)}(P)
\end{equation}
whenever $\al_k\in\ell^2$. Notice also that by Lemma \ref{L.1.4},
\begin{equation} \lb{3.1b}
\be^{(n)}(P)=\be(P)
\end{equation}
when $P\in M_m^\infty$ and $\omega(P)\ge m+n$. Since we have fixed
the $\al_k$'s, it will be more transparent to use the notation
$\be(P)$, $\be^{(n)}(P)$ rather than $\be(\mu,P)$,
$\be(\mu^{(n)},P-n)$.

Next we show that Theorem \ref{T.1.6} easily extends to
$D(\mu_n)$, $\Phi_n^*(\mu)$, and $\del D(\mu)$.

\begin{lemma} \lb{L.3.1}
For $m\ge 1$ and any $\mu$ we have
\begin{align}
w_{m}(\mu_n) & = \sum_{P\in M^n_m} N(P) \be(P) = w_{n+1,m}(\mu)
\lb{3.2}
\\ \del w_m(\mu) & =  \sum_{P\in M^\infty_m} N(P) \big[\be(P)-\be^{(1)}(P)\big]  \lb{3.3}
\end{align}
\end{lemma}

\begin{proof}
The first equality in \eqref{3.2} is nothing but \eqref{1.15} for
the measure $\mu_n$ instead of $\mu$. Then \eqref{1.3},
\eqref{1.1b}, and \eqref{1.1c} show that
\[
\log \Phi^*_{n+1}(\mu,z) = \log \Phi^*_{n+1}(\mu_n,z) =
\sum_{k=0}^{n} \log \rho_k - \log D(\mu_n,z)
\]
since $\Phi^*_{n+1}(\mu_n,z) =
\|\Phi^*_{n+1}(\mu_n,z)\|_{L^2(d\mu)} D^{-1}(\mu_n,z)$, and so
$w_{n+1,m}(\mu) = w_{m}(\mu_n)$ for $m\ge 1$.

By \eqref{3.0}, the $m$-th Taylor coefficient of $\del D(\mu,z)$
(and so of $\log \del D(\mu,z)$, too) only depends on $\al_0$, the
first $m$ Taylor coefficients of $f$ and first $m-1$ of $f_1$. That
is, $\del w_m(\mu)$ is a function of $\al_0,\dots,\al_m$ only (see
(1.3.48) in \cite{Sim1}). This means that for any $n\ge m$ we have
\[
\del w_m(\mu)=\del w_m(\mu_n) = w_m(\mu_n)-w_m((\mu_n)^{(1)}) =
\sum_{P\in M^n_m} N(P) \big[\be(P) -\be^{(1)}(P)\big]
\]
where the second equality is \eqref{3.0c} for $\mu_n$ and the
third follows from \eqref{3.2} and \eqref{3.1a}. But the last sum
equals the right hand side of \eqref{3.3} because \eqref{3.1b}
shows that $\be^{(1)}(P)=\be(P)$ when $P\in
M^\infty_m\smallsetminus M_m^n$.
\end{proof}

After this preparation we are ready to provide a characterization of
sequences of Verblunsky coefficients corresponding to measures $\mu$
for which $\log w(\tht)$ is integrable with respect to some
polynomial weight $|Q(e^{i\tht})|^2$. We let $Q(z) \equiv
\prod_{m=0}^N q_mz^m$ and define $p_m$ by
\[
|Q(z)|^2= \sum_{m=-N}^N p_m z^{m} = p_0 + \sum_{m=1}^N 2 \Real(p_m
z^{m}) \qquad \text{ for $|z|=1$}
\]
(note that $p_{m}= q_N\bar q_{N-m}+\dots+ q_m\bar q_0=\bar p_{-m}$).
With the convention $\log 0=-\infty$ we set
\begin{equation} \lb{3.4}
Z_Q(\mu)\equiv \int |Q(e^{i\tht})|^2 \log w(\tht)
\frac{d\tht}{4\pi}
\end{equation}
which is defined for any $\mu$ but can be $-\infty$. This is
because with $\log_\pm x\equiv\max\{\pm\log x, 0\}$,
\begin{equation} \lb{3.4a}
\int |Q(e^{i\tht})|^2 \log_+ w(\tht) \frac{d\tht}{4\pi} \le
\|Q\|^2_\infty \int w(\tht) \frac{d\tht}{4\pi} \le \frac
{\|Q\|^2_\infty}{4\pi}
\end{equation}
but the integral of $\log_- w(\tht)$ can be infinite, for instance,
when $w(\tht)=0$ on a set of positive measure. It is more common to
let $Z_Q$ be the negative of \eqref{3.4}, so that it is bounded from
below rather than above, but our definition will be more convenient
here. Note also that by \eqref{3.0z} with $p=1$ subsequently applied
to $\mu^{(n)}$, $n\ge 0$, in place of $\mu$, we have either
$Z_Q(\mu^{(n)})=-\infty$ for all $n\ge 0$ or
$Z_Q(\mu^{(n)})>-\infty$ for all $n\ge 0$.

Before determining the condition for $Z_Q(\mu)>-\infty$, we prove
the following {\it step-by-step sum rule}.

\begin{lemma} \lb{L.3.10}
For any $\mu$
\begin{equation} \lb{3.8a}
Z_Q(\mu) =  p_0 \log \rho_0 + \sum_{m=1}^N \Real \bigg( \bar p_m
\sum_{P\in M^\infty_m} N(P) \big[\be(P)-\be^{(1)}(P)\big] \bigg) +
Z_Q(\mu^{(1)})
\end{equation}
\end{lemma}

\begin{proof}
The sum on the right hand side of \eqref{3.8a} is always finite
since it has only finitely many non-zero elements and so
\eqref{3.8a} holds if both $Z_Q$ terms are $-\infty$. If both are
finite, then \eqref{3.0z} holds and so
\[
 \int |Q(e^{i\tht})|^2 \log \bigg( \frac{w(\tht)}{w^{(1)}(\tht)}
\bigg) \frac{d\tht}{4\pi} =  \frac 12 p_0 \del w_0(\mu) +
\sum_{m=1}^N \Real (\bar p_m \del w_m(\mu))
\]
together with \eqref{3.0e} and \eqref{3.3} gives \eqref{3.8a}.
\end{proof}

\begin{theorem} \lb{T.3.2}
For any $\mu$ and $Q$,
\begin{equation} \lb{3.5}
Z_Q(\mu) =  \sum_{k=0}^\infty   \Real \bigg( p_0\log \rho_k +
\sum_{m=1}^N \bar p_m \sum_{P\in M^0_m} N(P) \be(P+k) \bigg)
\end{equation}
%\begin{equation} \lb{3.5} Z_Q(\mu) =
%\lim_{n\to\infty} Z_Q(\mu_n)
%\end{equation}
%Moreover,
%\begin{equation} \lb{3.6}
%Z_Q(\mu_n) = - p_0 \sum_{k=0}^n \log \rho_k - \sum_{m=1}^N \Real
%\bigg( p_m \sum_{P\in M^n_m} N(P) \be(P) \bigg)
%\end{equation}
\end{theorem}

{\it Remark.} This shows that $Z_Q(\mu)$ is finite if and only if
the above sum converges.

\begin{proof}
Since \eqref{1.9b} and \eqref{3.4} give
\begin{equation} \lb{3.5x}
Z_Q(\mu_n) = \frac 12 p_0 w_0(\mu_n) + \sum_{m=1}^N \Real (\bar
p_m w_m(\mu_n))
\end{equation}
it follows from \eqref{1.9a} and \eqref{3.2} that
\begin{equation} \lb{3.6}
Z_Q(\mu_n) =  \sum_{k=0}^n  p_0\log \rho_k + \sum_{m=1}^N \Real
\bigg( \bar p_m \sum_{P\in M^n_m} N(P) \be(P) \bigg)
\end{equation}
Hence by \eqref{1.15} and $N(P)=N(P-\omega(P))$, the claim is
equivalent to $Z_Q(\mu) = \lim_{n\to\infty} Z_Q(\mu_n)$.

It is well known that $Z_Q$ is an entropy and therefore upper
semi-continuous in $\mu$ with respect to weak convergence of
measures (see Section 2.3 in \cite{Sim1}). In particular, since
$\mu_n\rightharpoonup\mu$, we obtain
\[
Z_Q(\mu) \ge  \limsup_{n\to\infty} Z_Q(\mu_n)
\]
Thus we are left with proving
\begin{equation} \lb{3.7}
Z_Q(\mu) \le  \liminf_{n\to\infty} Z_Q(\mu_n)
\end{equation}
This is obviously true if $Z_Q(\mu)= -\infty$, so assume that
$Z_Q(\mu^{(n)})>-\infty$ for all $n\ge0$. Then $w(\tht)>0$ for
a.e. $\tht$ and by Rakhmanov's theorem, $\al_n\to 0$ as
$n\to\infty$.

%We have
%\[
% \int |Q(e^{i\tht})|^2 \log \bigg( \frac{w(\tht)}{w^{(1)}(\tht)}
%\bigg) \frac{d\tht}{4\pi} =  \frac 12 p_0 \del w_0(\mu) +
%\sum_{m=1}^N \Real (\bar p_m \del w_m(\mu))
%\]
%which together with \eqref{3.0e} and \eqref{3.3} gives
%\[
% \int |Q(e^{i\tht})|^2 \log \bigg( \frac{w(\tht)}{w^{(1)}(\tht)}
%\bigg) \frac{d\tht}{4\pi} =  p_0 \log \rho_0 + \sum_{m=1}^N \Real
%\bigg( \bar p_m \sum_{P\in M^\infty_m} N(P)
%\big(\be(P)-\be^{(1)}(P)\big) \bigg)
%\]
%As mentioned before, the right hand side is always finite since
%the sum has only finitely many non-zero elements. By splitting the
%$\log$ in the left hand side integral we therefore obtain
%\begin{equation} \lb{3.8}
%Z_Q(\mu) =  p_0 \log \rho_0 + \sum_{m=1}^N \Real \bigg( \bar p_m
%\sum_{P\in M^\infty_m} N(P) \big(\be(P)-\be^{(1)}(P)\big) \bigg) +
%Z_Q(\mu^{(1)})
%\end{equation}
%where the $Z_Q$ terms are either both finite or both $-\infty$.

The step-by-step sum rule \eqref{3.8a} for $\mu^{(n)}$ in place of
$\mu$ reads
\[
Z_Q(\mu^{(n)}) = p_0 \log \rho_n + \sum_{m=1}^N \Real \bigg( \bar
p_m \sum_{P\in M^\infty_m} N(P)
\big[\be^{(n)}(P)-\be^{(n+1)}(P)\big] \bigg) + Z_Q(\mu^{(n+1)})
\]
and therefore we can iterate it and cancel the terms in the
telescoping sum to obtain
\[
Z_Q(\mu) = p_0 \sum_{k=0}^n \log \rho_k + \sum_{m=1}^N \Real \bigg(
\bar p_m \sum_{P\in M^\infty_m} N(P) \big[\be(P)-\be^{(n+1)}(P)\big]
\bigg) + Z_Q(\mu^{(n+1)})
\]
Using $\mu^{(n)}\rightharpoonup \tfrac{d\tht}{2\pi}$ (since
$\al_n\to 0$), $Z_Q(\tfrac{d\tht}{2\pi})=0$, and upper
semi-continuity of $Z_Q$, we obtain
\[
Z_Q(\mu) \le \liminf_{n\to\infty} \bigg[ p_0 \sum_{k=0}^n \log
\rho_k + \sum_{m=1}^N \Real \bigg( \bar p_m \sum_{P\in M^\infty_m}
N(P) \big[\be(P)-\be^{(n+1)}(P)\big] \bigg) \bigg]
\]

We claim that the quantity inside the $\liminf$ differs by $o(1)$
from \eqref{3.6} (in which case \eqref{3.7} holds and we are
done). Indeed --- the difference of these two is at most
\[
\sum_{m=1}^N |p_m| \bigg( \sum_{P\in M^{n+m}_m\smallsetminus M_m^n}
|N(P)| |\be(P) - \be^{(n+1)}(P)| \bigg)
\]
by \eqref{3.1b} and the fact that $\be^{(n+1)}(P)=0$ when $P\in
M_m^n$. This sum has a uniformly bounded number of terms for all
$n$, both $p_m$ and $N(P)$ are also bounded by a constant not
depending on $n$ (only on $N$ and $Q$), and
\[
\lim_{n\to\infty} \sup_{P\in M^{n+m}_m\smallsetminus M_m^n} \{
|\be(P)|+ |\be^{(n+1)}(P)| \} = 0
\]
since $\al_n\to 0$.
\end{proof}

We have thus expressed $Z_Q(\mu)$ as an infinite sum in terms of
the Verblunsky coefficients of $\mu$. We can now apply Theorem
\ref{T.1.6x} to prove Theorem \ref{T.3.3}.

\begin{proof}[Proof of Theorem \ref{T.3.3}]
The right hand side of \eqref{3.10} is equivalent to
$Z_Q(\mu)>-\infty$. By Theorem~\ref{T.3.2}, this happens precisely
when $\lim_{n\to\infty}$ \eqref{3.5x} $>-\infty$. But
\[
\frac 12 p_0w_0(\mu_n)=p_0\sum_{k=0}^n \log\rho_k = -\frac 12
p_0\sum_{k=0}^n \big( |\al_k|^2 + O(|\al_k|^4) \big)
\]
and by Theorem \ref{T.1.6x} applied to $\mu_n$, the sum in
\eqref{3.5x} is equal to
\[
\sum_{m=1}^N \Real \bigg( \bar p_m \bigg( R_m(\mu_n)-\sum_{k\le
n-m} \al_{k+m} \bar\al_{k} \bigg)\bigg)
\]
%\begin{equation} \nonumber %\lb{3.10}
%\[
%Z_Q(\mu_n)=-\bigg[ \frac 12 p_0\sum_{k\le n} |\al_k|^2 +
%\sum_{m=1}^N \Real \bigg( \bar p_m \sum_{k\le n-m} \al_{k+m}
%\bar\al_{k} \bigg) \bigg] + \sum_{m=1}^N \Real (\bar p_m
%R_m(\mu_n)) + O(1)
%\]
%\end{equation}
%\[
%\sum_{m=1}^N \Real \bigg( \bar p_m \bigg( R_m(\mu_n)-\sum_{k\le
%n-m} \al_{k+m} \bar\al_{k} \bigg)\bigg)
%\]
The estimate \eqref{1.17} and the hypothesis show that $R_m(\mu_n)$
and $\sum_{k=0}^n O(|\al_k|^4)$ are uniformly bounded in $n$, so it
only remains to show that $\{\bar Q(S)\al\}_k\in\ell^2$ is
equivalent to
\begin{equation} \lb{3.11}
-\bigg[ \frac 12 p_0\sum_{k\le n} |\al_k|^2 + \sum_{m=1}^N \Real
\bigg( \bar p_m \sum_{k\le n-m} \al_{k+m} \bar\al_{k} \bigg)
\bigg]
\end{equation}
being uniformly bounded in $n$. We write
\begin{align*}
\sum_{k\le n} & |\{\bar Q(S)\al\}_k|^2  = \sum_{k\le n} |\bar q_N
\al_{k+N} +\dots+ \bar q_0 \al_k|^2
\\ & = (|q_N|^2+\dots+|q_0|^2)\sum_{k\le n} |\al_k|^2 +
\sum_{m=1}^N 2\Real \bigg( (\bar q_Nq_{N-m}+\dots+\bar q_mq_0)
\sum_{k\le n} \al_{k+m}\bar\al_k \bigg) + O(1)
\\ & = p_0 \sum_{k\le n} |\al_k|^2 + \sum_{m=1}^N 2\Real \bigg( \bar p_m
\sum_{k\le n-m} \al_{k+m}\bar\al_{k} \bigg) + O(1)
\end{align*}
In the second equality the remainder $O(1)$ is bounded by a
constant independent of $n$ because it is a sum of a bounded
number of terms involving only $\al_k$ with $k\le N$ or $|k-n|\le
N$. And the last equality holds because $p_m = q_N\bar
q_{N-m}+\dots+ q_m\bar q_0$ and $\sum_{k= n-m+1}^n
\al_{k+m}\bar\al_{k}$ is uniformly bounded in $n$ and so $O(1)$.
Hence $\{\bar Q(S)\al\}_k\in\ell^2$ if and only if \eqref{3.11} is
uniformly bounded in $n$.
\end{proof}

Recall that in the proof of Theorem \ref{T.3.2} we have showed
$Z_Q(\mu) = \lim_{n\to\infty} Z_Q(\mu_n)$. Here is a
generalization of this fact.

\begin{proposition} \lb{C.3.4}
If $f=gQ$ with $g\in C(\partial\mathbb D)$ a positive function and
$Q$ a polynomial, then for any $\mu$,
\begin{equation} \lb{3.12}
Z_f(\mu) = \lim_{n\to\infty} Z_f(\mu_n)
\end{equation}
\end{proposition}

\begin{proof}
We again have $Z_f(\mu)\ge\limsup_{n} Z_f(\mu_n)$ because $Z_f$ is
upper semi-continuous as well \cite{Sim1}. Let $g_\eps$ be a
polynomial such that $g^2\le |g_\eps|^2\le g^2+\eps$ on
$\partial\bbD$. Such a polynomial exists because the functions
$h(z)=\sum_{k=-K}^K c_kz^k$ are dense in $C(\partial\bbD)$ (by the
complex Stone-Weierstrass theorem) and the polynomial $z^K h(z)$
satisfies $|z^K h(z)|=|h(z)|$ for $|z|=1$. Then $\lim_{\eps\to 0}
Z_{g_\eps Q}(\mu)= Z_f(\mu)$ because $g$ is bounded away from $0$,
and
\[
Z_{f}(\mu_n) \ge Z_{g_\eps Q}(\mu_n) - \frac \eps{4\pi}
\|Q\|^2_\infty
\]
by \eqref{3.4a}. Since $g_\eps Q$ is a polynomial, the proof of
Theorem \ref{T.3.2} shows $\lim_n Z_{g_\eps Q}(\mu_n)=Z_{g_\eps
Q}(\mu)$, and we obtain $\liminf_{n} Z_f(\mu_n)\ge Z_f(\mu)$ by
taking $\eps\to 0$.
\end{proof}

It is an interesting open question whether this result holds for
any $f$, not just such that vanish at only finitely many points of
$\partial\bbD$ and to an even degree.

%%%%%%%%%%%%%%%%%%%%%%%%%%%%%%%%%%%%%%%%%%%%%%%%%%%%%
\section{Relative Ratio Asymptotics} \lb{S4}
%%%%%%%%%%%%%%%%%%%%%%%%%%%%%%%%%%%%%%%%%%%%%%%%%%%%%

In this section we provide another application of our methods. We
prove Theorem \ref{T.4.1} and give a simple proof of a deep result,
in part due to Khrushchev \cite{Khr} and in part to Barrios and L\'
opez \cite{BL}, on ratio asymptotics $\Phi_{n+1}^*/\Phi_n^*$ as
$n\to\infty$ of the reversed polynomials (see also \cite[Section
9.5]{Sim2}). We also give a generalization of this result.
% and we derive a new result on relative ratio
%asymptotics $(\Phi_{n+1}^*/\Phi_n^*)(\Psi_{n+1}^*/\Psi_n^*)^{-1}$
%where $\Phi_n^*$ and $\Psi_n^*$ are the reversed polynomials of
%measures $\mu$ and $\nu$, respectively.

\begin{proof}[Proof of Theorem \ref{T.4.1}]
Let us define
\[
\Omega_n(\mu,\nu)\equiv \frac{\Phi_{n+1}^*(\mu) / \Phi_{n}^*(\mu)}
{\Phi_{n+1}^*(\nu) / \Phi_{n}^*(\nu)}, \qquad
\log\Omega_n(\mu,\nu)\equiv\sum_{m\ge1} \omega_{n,m}(\mu,\nu)z^m
\]
(recall that $\Phi_n^*(0)=1$). It follows from
$|\Phi_n(z)|\le|\Phi_n^*(z)|$ for $z\in\bbD$ (see (1.7.1) in
\cite{Sim1}) and from
\begin{equation} \lb{4.1a}
\frac{\Phi_{n+1}^*}{\Phi_n^*} = 1-\al_n z \frac{\Phi_n}{\Phi_n^*}
\end{equation}
(see \eqref{1.3}) that the ratio $\Phi_{n+1}^*/\Phi_n^*$ is bounded
away from $0$ and $\infty$ on any compact $K\subset\bbD$. Hence
\eqref{4.1} is equivalent to $\Omega_n(\mu,\nu)\to 1$ as
$n\to\infty$ uniformly on compact subsets of $\bbD$, which in turn
is equivalent to $\omega_{n,m}\to 0$ as $n\to\infty$ for each
$m\ge1$.

%We first note that both the numerator and denominator in
%\eqref{4.1} are bounded away from $0$ and $\infty$ on any compact
%$K\subset\bbD$. This follows from $|\Phi_n(z)|\le|\Phi_n^*(z)|$
%for $z\in\bbD$ (see (1.7.1) in \cite{Sim1}) and from
%\begin{equation} \lb{4.1a}
%\frac{\Phi_{n+1}^*}{\Phi_n^*} = 1-\al_n z \frac{\Phi_n}{\Phi_n^*}
%\end{equation}
%(see \eqref{1.3}). Thus the convergence to 1 of \eqref{4.1} on
%compacts inside $\bbD$ is equivalent to that of $\log \eqref{4.1}$
%to 0, which in turn is equivalent to the convergence of Taylor
%coefficients of $\log \eqref{4.1}$ to 0 (because $\log$ \eqref{4.1}
%is also uniformly bounded on $K$).
%We will be done if we prove that the convergence of
%the first $m$ Taylor coefficients of $\log\eqref{4.1}$ is equivalent
%to that of $\al_n\bar\al_{n-\ell} - \ga_n\bar\ga_{n-\ell}$ as
%$n\to\infty$ for all $\ell\le m$.

Now
\[
\log \Omega_n(\mu,\nu) = \log \Phi_{n+1}^*(\mu) - \log
\Phi_{n}^*(\mu) - \log \Phi_{n+1}^*(\nu) + \log \Phi_{n}^*(\nu)
\]
and so by \eqref{3.2} and \eqref{1.15a}
\begin{equation} \lb{4.3}
\omega_{n,m}(\mu,\nu)=\sum_{P\in M^0_m} N(P) \big( \be(\nu,P+n) -
\be(\mu,P+n) \big)
\end{equation}
for $m\ge 1$. If \eqref{4.2} holds, then obviously $\be(\nu,P+n) -
\be(\mu,P+n)\to 0$ for each $P\in M^0_m$, that is,
$\omega_{n,m}(\mu,\nu)\to 0$ as $n\to\infty$.

Assume now that $\omega_{n,m}(\mu,\nu)\to 0$ as $n\to\infty$. When
$m=1$, then $M^0_m=\{\{(0,1)\}\}$ and \eqref{4.3} equals just
$\al_n(\mu)\bar\al_{n-1}(\mu) - \al_n(\nu)\bar\al_{n-1}(\nu)$. Hence
\eqref{4.2} holds for $\ell=1$. We proceed by induction, so assume
\eqref{4.2} holds for $\ell=1,\dots,m-1$. If $P\in M^0_m$ and
$|P|\ge2$, then by the induction hypothesis $\be(\nu,P+n) -
\be(\mu,P+n)\to 0$ (because $\max\{ a_l\,|\, (k_l,a_l)\in P\}\le m
-1$). Since \eqref{4.3} converges to 0 and the only element of
$M^0_m$ with $|P|=1$ is $P=\{(0,m)\}$ (in which case $N(P)=-1$), it
follows that $\be(\nu,\{(n,m)\}) - \be(\mu,\{(n,m)\})\to 0$ as well.
But this is \eqref{4.2} for $\ell=m$.
\end{proof}

For any $a\in [0,1]$ we define
\begin{equation*} %\lb{4.3a}
G_a(z)\equiv \frac 12 \,\big( 1+z+\sqrt{(1-z)^2+4a^2z}\, \big)
\end{equation*}
with the usual branch of the square root (in particular,
$G_0\equiv 1$). We then have

\begin{theorem}[\cite{Khr} and \cite{BL}] \lb{T.4.2}
Let $\mu$ be a non-trivial probability measure on $\partial\bbD$.
Then
\begin{equation} \lb{4.4}
\Omega_n(\mu)\equiv \frac {\Phi_{n+1}^*} {\Phi_{n}^*}
\end{equation}
converges uniformly on compact subsets of $\bbD$ as $n\to\infty$ if
and only if for each $\ell\ge 1$ there is $c_\ell\in\overline{\bbD}$
such that
%there is $c_\ell$ such that
\begin{equation} \lb{4.5}
\lim_{n\to\infty}  \al_n\bar\al_{n-\ell}  = c_\ell
\end{equation}
Moreover, if \eqref{4.5} holds for all $\ell\ge 1$, then
%either all $c_\ell=0$ and the limit of \eqref{4.4} is $1$, or
$c_\ell = a^2\la^\ell$ for some $a\in [0,1]$ and $|\lambda|=1$ and
\[
\lim_{n\to\infty}\Omega_n(\mu;z)=G_a(\lambda z)
\]
\end{theorem}

{\it Remark.} In particular, $\lim_n\Omega_n(\mu)=1$ precisely
when all $c_\ell=0$. In this case \eqref{4.5} is called {\it M\'
at\' e-Nevai condition}. Accordingly, one might call \eqref{4.2}
{\it relative M\' at\' e-Nevai condition}.

\begin{proof}
Equivalence of the convergence of \eqref{4.4} and \eqref{4.5} is
proved in the same way as Theorem \ref{T.4.1}. The only difference
is that now with
\[
\log\Omega_n(\mu)\equiv\sum_{m\ge1} \omega_{n,m}(\mu)z^m
\]
\eqref{4.3} reads
\begin{equation} \lb{4.5a}
\omega_{n,m}(\mu)=-\sum_{P\in M^0_m} N(P) \be(\mu,P+n)
\end{equation}
and ``$\be(\nu,P+n) - \be(\mu,P+n)\to 0$'' and
``$\omega_{n,m}(\mu,\nu)\to 0$'' are replaced by ``$\be(\mu,P+n)$
converges'' and ``$\omega_{n,m}(\mu)$ converges'', respectively, in
the argument (we actually have $\be(\mu,P+n)\to c_{a_1}\dots
c_{a_i}$ when $P=\{(k_1,a_1),\dots,(k_i,a_i) \}$).
%and therefore the respective limits need not only be 1 in
%\eqref{4.4} and 0 in \eqref{4.5}.
Note that the proof also shows that $\lim_n\Omega_n(\mu)=1\equiv
G_0$ precisely when all $c_\ell=0$ (and so $a=0$).

%Clearly, \eqref{4.5} implies
%\begin{equation} \lb{4.5x}
%\lim_{n\to\infty} \omega_{n,m}(\mu)=-\sum_{P\in M^0_m} N(P)
%c_{a_1}c_{a_2}\ldots c_{a_i}, \qquad P=\{(k_1,a_1),\dots,(k_i,a_i)
%\}
%\end{equation}
%and so $\Omega_n(\mu)$ converges. Conversely, the convergence of
%$\omega_{n,m}$ as $n\to\infty$ for each $m\ge1$ implies
%\eqref{4.5} by induction as above. Moreover, since
%$\omega_{n,1}(\mu)=\al_n\bar\al_{n-1}$, then $\omega_{n,1}\to
%c_1=0$ yields $c_1=c_2=\ldots=0$ (cf. Lemma 2.2 in \cite{Khr}).
%Obviously, $c_1=c_2=\ldots=0$ gives $\omega_{n,m}(\mu)\to 0$ as
%$n\to\infty$ for each $m\ge1$ by \eqref{4.5x} and so $\Omega_n\to
%1$. The structure formula for $c_\ell$'s is proved in \cite [Lemma
%2.2] {Khr}, so it
%so that we use $\be(\mu,P+n)\to c_{a_1}\dots c_{a_i}$ as
%$n\to\infty$ when $P=\{(k_1,a_1),\dots,(k_i,a_i) \}$, and
%therefore the respective limits need not only be 1 in \eqref{4.4}
%and 0 in \eqref{4.5}. However, the proof also shows that the limit
%of \eqref{4.4} is $G_0\equiv 1$ precisely when all $c_\ell=0$ (and
%so $a=0$).

Hence assume \eqref{4.5} holds with not all $c_\ell=0$. It is
obvious that if $c_1=0$, then the existence of the limit $c_\ell$
implies $c_\ell=0$ for all $\ell$. Thus we must have
$c_1=a^2\lambda$ for some $a\in(0,1]$ and $|\lambda|=1$. In
particular, $\liminf_n |\al_n |>0$. But then $|\al_{n+3}||\al_{n+2}|
- |\al_{n+1}||\al_{n}|\to 0$ and $|\al_{n+3}||\al_{n+1}| -
|\al_{n+2}||\al_{n}|\to 0$ (both by \eqref{4.5}) give $|\al_{n+2}| -
|\al_{n+1}|\to 0$, which together with $|\al_{n+2}||\al_{n+1}|\to
a^2$ gives $|\al_{n}|\to a$. This and \eqref{4.5} imply
$\al_{n+1}\al_n^{-1}\to\la$, and then
$\al_n\al_{n-\ell}^{-1}\to\la^\ell$ so that $c_\ell=a^2\lambda^\ell$
for all $\ell$.

It remains to prove that in the case $a\not=0$ the limit $G(z)$ of
\eqref{4.4} is $G_a(\la z)$. Let $\nu$ be the measure with
Verblunsky coefficients $\al_n(\nu)\equiv a\la^n$ if $0<a<1$ and
$\al_n(\nu)\equiv a_n\la^n$ with $a_n\uparrow 1$ if $a=1$. Then
Theorem~\ref{T.4.1} applies and so the limit function of
$\Phi_{n+1}^*(\nu)/\Phi_n^*(\nu)$ is also $G(z)$. By \eqref{4.1a}
we know that the limit $H(z)$ of
$\al_n(\nu)\Phi_n(\nu)/\Phi_n^*(\nu)$ also must exist and
\begin{equation} \lb{4.6}
G=1-zH
\end{equation}
From \eqref{1.2} we have
\[
\la z\al_n(\nu)\,\frac{\Phi_n(\nu)}{\Phi_n^*(\nu)} =
\la\al_n(\nu)\,\frac{\Phi_{n+1}(\nu)}{\Phi_n^*(\nu)}+a^2\la =
\al_{n+1}(\nu)\,\frac{\Phi_{n+1}(\nu)}{\Phi_{n+1}^*(\nu)}
\frac{\Phi^*_{n+1}(\nu)}{\Phi_n^*(\nu)} + a^2\la
\]
Therefore $\la zH= HG+a^2\la$. We substitute $H(G-\la z)=-a^2\la$
into \eqref{4.6} multiplied by $G-\la z$ to obtain $G^2-(1+\la z)G
+(1-a^2)\la z=0$. Using $G(0)=1$, it follows that $G(z)=G_a(\la
z)$.
\end{proof}

We conclude with using \eqref{4.5a} to obtain the following
generalization of Theorem \ref{T.4.2}.

\begin{theorem} \lb{T.4.3}
Let $\mu$ be a non-trivial probability measure on $\partial\bbD$,
let $\{j_n\}$ be an increasing sequence of integers and
$k'\in\bbZ\cup\{\infty\}$. Then $\Omega_{k+j_n}(\mu)$ from
\eqref{4.4} converges for any $k<k'$ uniformly on compact subsets of
$\bbD$ as $n\to\infty$ if and only if for each $\ell\ge 1$ and
$k<k'$ there is $c_{k,\ell}\in\overline{\bbD}$ such that
\begin{equation} \lb{4.7}
\lim_{n\to\infty}  \al_{k+j_n}\bar\al_{k+j_n-\ell}  = c_{k,\ell}
\end{equation}
\end{theorem}

{\it Remark.} For the special case $j_n=np$ with $p\ge 1$ see
\cite[Theorem 9.5.10]{Sim2}.

\begin{proof}
We again follow the lines of the two previous proofs. In one
direction we have that the existence of all the $c_{k,\ell}$ with
$k<k'$ gives the convergence of $\be(\mu,P+k+j_n)$ for any $P\in
\bigcup_mM_m^0$ and $k<k'$ as $n\to\infty$ (note that if
$(k_l,a_l)\in P\in M_m^0$, then $k_l\le 0$). This in turn gives the
convergence of $\omega_{k+j_n,m}(\mu)$ for any $k<k'$ and $m$ by
\eqref{4.5a}, and thus that of $\Omega_{k+j_n}(\mu)$ for $k<k'$.

In the opposite direction, the reverse of this argument and
induction on $\ell$ as in the proof of Theorem \ref{T.4.1} shows
that the convergence of $\Omega_{k+j_n}(\mu)$ for any $k<k'$ implies
\eqref{4.7}. This again uses the fact that if $(k_l,a_l)\in P\in
M_m^0$, then $k_l\le 0$.
\end{proof}

\bigskip
%%%%%%%%%%%%%%%%%%%%%%%%%%%%%

\end{document}